\documentclass[12pt]{article}
   \usepackage{amscd,amsthm,amssymb,amsfonts,amsmath}

\theoremstyle{plain}

\newtheorem{lemma}{Lemma}
\newtheorem{prop}{Proposition}
\newtheorem{cor}{Corollary}

\theoremstyle{definition}
\newtheorem{defn}{Definition}

\theoremstyle{remark}
\newtheorem*{remark}{Remark}

\newcommand{\CC}{\mathbb C}    
 
\newcommand{\gl}{\mathfrak {gl}}    

\makeatletter
\renewcommand{\@makefnmark}{\mbox{\textsuperscript{}}}

\begin{document}

\title{Symmetric functions for the  generating matrix of  Yangian of $\gl_n(\CC)$.}
\maketitle
\begin{center}
Natasha Rozhkovskaya\\

Department of Mathematics\\ 
Kansas State University, USA\\

\end{center}

\begin{abstract} Analogues of classical combinatorial identities for  elementary and homogeneous symmetric functions with  coefficients in Yanigian are discussed. As a corollary, similar relations are deduced for shifted Schur functions.
\end{abstract}
\section*{Introduction}
In this note we  prove some combinatorial relations between the analogues of  symmetric functions for the Yangian   of the Lie algebra $\gl_n(\CC)$. 
The applications of the results are illustrated by deducing  properties of   Capelli polynomials and shifted symmetric functions. Some of these properties  were obtained, for example,  in \cite{OO1} from  the definitions of shifted symmetric functions. Here, due to  the existence of  evaluation homomorphism, they become immediate consequences of similar combinatorial formulas in the Yangian. 
The elementary symmetric functions in the Yangian   of the Lie algebra $\gl_n(\CC)$ are known to be generators of Bethe sublagebra. Bethe subalgebra  finds numerous applications in quantum integrable models of XXX type and Gaudin type (\cite{MTV1}, \cite{MTV2}, \cite{MTV3}). We  describe  the inverse of the universal differential operator for higher transfer matrices of XXX model.

The author is very grateful  to E.\,Mukhin for encouraging discussions and valuable remarks. She would  like to thank
 P.\,Pyatov,  A.\,Chervov for sharing comments on the subject. The hospitality   of Institut des Hautes Etudes Scientifiques  and  of Max Plank Institute  for Mathematics in Bonn provided very inspiring 
 atmosphere for the research.
 The project is supported in part by KSU Mentoring fellowship for WMSE.
\section*{Notations and Preliminary facts}
The following  notations will be used through the paper.
All non-commutative   determinants are defined to be row determinants. Namely, if $X$ is a matrix with entries $(x_{ij})_{i,j=1,\dots n}$ in an  associative algebra $A$, put
$$
\text{det} X=\text{rdet} X=\sum_{\sigma\in S_n}(-1)^{\sigma}x_{1\sigma(1)}\dots x_{n\sigma(n)},
$$
where the sum is taken over all permutations of $n$ elements. We also define the following  types  of powers of the matrix $X$:

$$X^{[k]}:= X_1\dots X_k\, \in\, \text {End}(\CC^n)^{\otimes k}\otimes A,$$
where 
$$
X_s=\sum_{ij}1\otimes\dots\otimes \underset{s}{E_{ij}}\otimes\dots \otimes 1\otimes x_{ij},
$$
and $$X^{k}:= X\dots X\, \in\, \text {End}(\CC^n)\otimes A.$$
(This is  just regular multiplication of matrices).

\subsection*{Definition of Yangian}

Let $P_{l,m}$ be a permutation of $l$-th and $k$-th copies of $\CC^n$ in $(\CC^n)^{\otimes k}$:
\begin{equation}\label{perm}
P_{l,m}= \sum_{ij}1\otimes \dots \otimes1\otimes\underset{l}{ E_{ij}}\otimes \dots 
\otimes \underset{m}{E_{ji}}\otimes 1\dots \otimes 1.
\end{equation}
Let $u$ be an independent variable. Consider the Yang matrix
$$
R(u)=1-\frac{P_{12}}{u}\quad  \in \text{End}\,(\CC^n)^{\otimes 2}[[u^{-1}]].
$$
\begin{defn}The Yangian $Y(n)$ of the Lie algebra $\gl_n(\CC)$ is an associative unital algebra, generated by the elements $\{ t_{ij}^{(k)}\}$,
($i,j=1\dots n$, $k=1,2,...$), satisfying  the  relation
\begin{equation}\label{RTT}
R(u-v)T_1(u)T_2(v)= T_2(v) T_1(u)R(u-v).
\end{equation} 
Here $T(u)=(t_{ij}(u))_{1\le i,j\le k}$  is the generating matrix of  $Y(n)$: the entries of  $T(u)$ are formal power series with coefficients in $Y(n)$:
$$
 t_{ij}(u)=\sum_{k=0}^{\infty} \frac{t_{ij}^{(k)}}{u^{k}},\quad t_{ij}^{(k)}\in Y(n),  \quad t_{ij}^{(0)}=\delta_{i,j}.
$$
\end{defn}
 The definition of $Y(n)$ implies that many formulas involving its generating matrix $T(u)$  contain the shifts of the parameter $u$.
To simplify  some of these formulas, it is convenient to  introduce a shift-variable $\tau$ (we follow \cite{Tal},\cite{ MTV1}, \cite{Cherv}  in this approach). 
Any element $f(u)$ of $Y(n)[[u^{-1}]]$ we identify with the operator of multiplication by this formal power series, acting on $Y(n)[[u^{-1}]]$.
Let $\tau^{\pm 1}=e^{\frac{\pm d}{du}}$. These operators also act on $Y(n)[[u^{-1}]]$ by shifts of the variable $u$:
\begin{equation}\label{Euler}
\tau^{\pm}(g(u))= e^{\frac{\pm d}{du}}(g(u))=g(u\pm1), \quad g(u)\in Y(n)[[u^{-1}]].
 \end{equation}
Thus, under this identification of shifts
$\tau^{\pm 1}$   and the elements $f(u)$ of $Y(n)[[u^{-1}]]$ with differential operators acting on the algebra
$Y(n)[[u^{-1}]]$, we can write the following commutation relation:
\begin{equation}\label{tau}
\tau^{\pm}\,f(u)=f(u\pm1)\,\tau^{\pm}.
\end{equation}
We will use the relation (\ref{tau}) to write the formulas for symmetric  functions 
$e_k(u,\tau)$, $h_k(u,\tau)$, $p_k(u,\tau)$, defined in the next section.

\subsection*{Symmetrizer and antisymmetrizer.}  Define the  projections to the symmetric and antisymmetric part of $(\CC^n)^{\otimes m}$:
 $$A_k=\frac{1}{k!}\sum_{\sigma\in \mathcal{S}_k}(-1)^\sigma \sigma,\quad \quad S_k=\frac{1}{k!}\sum_{\sigma\in \mathcal{S} _k} \sigma.$$  These are the elements of  the group algebra $\CC[\mathcal{S}_k]$ of the permutation group,  acting on $(\CC^n)^{\otimes k}$ by permuting the tensor components.
The operators  enjoy  the listed below properties.
\begin{prop}\label{AandS}
(a) $$A_k^2= A_k\quad \quad \text{  and  }\quad \quad S_k^2=S_k.$$

(b) With abbreviated notations   $R_{ij}=R_{ij}(v_i-v_j)$, write 
$$R(v_1,\dots v_m)=(R_{m-1,m})(R_{m-2,m}R_{m-2,m-1})...(R_{1,m}
\dots R_{1,2}).
$$
Then $A_k=\frac{1}{k!}R(u,u-1,\dots u-k+1)$, and 
$S_k=\frac{1}{k!}R(u,u+1,\dots u+k-1)$.

(c) $$A_kT_1(u)\dots T_m(u-k+1)=T_k(u-k+1)\dots T_1(u) A_k,$$
 $$S_k T_1(u)\dots T_k(u+k-1)=T_k(u+k-1)\dots T_1(u)S_k.$$
 
 (d) $$\text{tr}\,(A_nT_1(u)\dots T_n(u-n+1))=\text{qdet}\, T(u).$$

(e) $$A_{k+1}=\frac{1}{k+1}A_{k}
 \, R_{k,k+1}\left(\frac{1}{k}\right)\, A_k,$$
$$S_{k+1}=\frac{1}{k+1}S_{k}\,R_{k,k+1}\left(-\frac{1}{k}\right)\,S_k.$$

(h) Put $$B^{\mp}_l:=\frac{1}{l!} \, R_{l-1,l}\left(\frac{\pm1}{l-1}\right)R_{l-2,l-1}\left(\frac{\pm 1}{l-2}\right)\dots R_{1,2}\left(\pm1\right).$$
Then 
$$
S_k=B^+_{2}B^+_{3}\dots B^+_{k},\quad
A_k=B^-_{2}B^-_{3}\dots B^-_{k},
$$
\end{prop}

\begin{proof} The  properties (a) -- (d) are contained in Propositions 2.9 -- 2.11 in  \cite{Mol}. The property (e) can  be shown by induction.  The statement of (h) follows from  (e). Note that  (b) and (h) give different presentations of symmetrizer and antisymmetrizer in terms of R-matrices.
For example, by (b), $A_3=\frac{1}{6}R_{23}(1)\,R_{13}(2)\,R_{12}(1)$, and by property (h), 
$A_3=\frac{1}{12}R_{12}(1)\,R_{23}\left(\frac{1}{2}\right)\,R_{12}(1)$.
\end{proof}

\section*{Elementary and homogeneous symmetric functions}
\begin{defn}
The following   formal power sums in $u^{-1}$ with coefficients in $Y(n)$ are  the analogues of ordinary symmetric functions:

 {\it Elementary symmetric functions: }$$e_k(u)=tr\,( A_k T_1(u)\dots T_k(u-k+1)),$$

{\it Homogeneous symmetric functions:}  $$h_k(u)=tr\,( S_k T_1(u)\dots T_k(u+k-1)),$$

{\it Power sums:}
$$
 {p}^{\pm}_{ k}(u)=\text{tr}\left(T(u) T(u\pm1)\dots T(u\pm (k-1)\right)
.$$
\end{defn}
 
\subsection*{Bethe subalgebra}
Let $Z$ be a matrix of size $n$ by $n$ with complex coefficients.
Consider  $\mathcal{B}(\gl_n(\CC,Z))$  -- the  
 commutative subalgebra  of  the Yangian $Y(n)$, generated by the coefficients of all the series
 $$
b_k(u,Z)=\text{tr}\,(A_nT_1(u)\dots T_k(u-k+1)Z_{k+1}\dots Z_n),\quad k=1,2\dots n.
$$
It is called {\it Bethe subalgebra} 
(see, for example \cite{KBI}, \cite{KR}, \cite {KS}, \cite{Bet}).
The introduced above 
elements $e_k(u)$ are proportional to generators of the Bethe subalgebra  with    $Z$ being  the identity matrix:
\begin{lemma}$e_{k}(u)=\frac{n!}{k!\,(n-1)^{n-k}}\, b_k(u,\text{Id})$
\end{lemma} 
\begin{proof} 
Let $\text{tr}_{(1\dots a)}$ denote the trace by the first $a$ components in  the tensor product $(\text {End}\, (\CC^n))^{\otimes (m+1)}$.
By  Proposition \ref{AandS}  (c), (e), and the cyclic property of the trace, we obtain that 
\begin{align} \label{tr1}
\text{tr}_{(1\dots m+1)}\left(A_{m+1}\,T_1(u)\dots T_k(u-k+1)\otimes 1^{\otimes ^{m+1-k}}\right)\\\notag
={\frac{(n-1)}{m+1}} \,\text{tr}_{(1\dots m)}
 \left(A_m\,T_1(u)\dots T_k(u-k+1)\otimes 1^{\otimes ^{m-k}}\right).
\end{align}
From (\ref{tr1}) one can show by induction that 
\begin{align} \notag
b_k (u,\text{Id})=
\text{tr}_{(1\dots n)}
(A_{n}T_1(u)\dots T_k(u-k+1)\otimes 1^{\otimes ^{n-k}})
=\frac{(n-1)^{n-k}k!}{n!}\,e_k(u).
\end{align}
\end{proof}
\begin{prop}\label{B}
 Let the matrices $B^{\pm}_k$  be defined as in Proposition\,\ref{AandS}, (h).
Then
\begin{equation}\label{e_B}\begin{aligned}
e_k(u) =&\text{tr}\, \left(B^-_k\, T_1(u)\dots T_k(u-k+1)\right),\\
h_k(u)=&\text{tr}\, \left(B^+_k\,T_1(u)\dots T_k(u+k-1)\right),\\
e_k(u+k-1)&=\text{tr}\, \left(A_k\, T_1(u)\dots T_k(u+k-1)\right),\\
h_k(u-k+1)&=\text{tr}\, \left(S_k\,T_1(u)\dots T_k(u-k+1)\right).
\end{aligned}\end{equation}

\end{prop}
\begin{proof}
 By  Proposition \ref{AandS} part (e), 
\begin{equation}\label{proof2-4}
\begin{aligned}
e_k(u)&=\frac{1}{k}\text{tr}\,\left(A_{k-1}
\, R_{k-1,k}\left(\frac{1}{k-1}\right)\, A_{k-1}T_1(u)\dots T_k(u-k+1)\right),\\
 \quad &=\frac{1}{k}\text{tr}\,\left(\,R_{k-1,k}\left(\frac{1}{k-1}\right)\,A_{k-1}
T_1(u)\dots T_k(u-k+1)\,A_{k-1}\,  \right),\\
\quad & =\frac{1}{k}\text{tr}\,\left(\,R_{k-1,k}\left(\frac{1}{k-1}\right)\,\,A_{k-1}
 T_1(u)\dots T_k(u-k+1)  \right).\\
 \end{aligned}
\end{equation}
The last equality follows from properties (c) and (a) of the Proposition \ref{AandS}.
Applying  the same Proposition \ref{AandS} part (e) to $A_{k-1}$,
 and observing, that $A_{k-2}$ commutes with $R_{k-1,k}\left(\frac{1}{k-1}\right)$, 
we obtain that 
\begin{equation*}
\begin{split}
e_k(u)=\frac{1}{k(k-1)}\text{tr}\left(R_{k-1,k}\left(\frac{1}{k-1}\right) R_{k-2,k-1}\left(\frac{1}{k-2}\right) A_{k-2}
 T_1(u)\dots T_k(u-k+1)  \right).\\
 \end{split}
\end{equation*}
Proceeding by induction, we obtain the first  statement of (\ref{e_B}).
The second formula is proved similarly, and the last two  can be checked directly.
\end{proof}
Introduce  the following notations:
   \begin{equation}
   \begin{aligned}
       e_k(u,\tau)&=\text{tr}\left(A_k(T(u)\tau^{-1})^{[k]}\right),\\
       { h}_{k}(u,\tau)&=\text{tr}\left((S_kT(u) \tau)^{[k]}\right),\\
      {p}^{\pm}_{ k}(u,\tau)&=\text{tr}\left((T(u) \tau^{\pm 1})^{k}\right).
      \end{aligned}
\end{equation}
Observe that 
\begin{equation}
\begin{split}
e_k(u,\tau)=e_k(u)\tau^{-k},\quad 
h_k(u,\tau)=h_k(u)\tau^{k},\quad 
{p}^{\pm}_{ k}(u,\tau)={p}^{\pm}_ k(u) \tau^{\pm k}.
\end{split}
\end{equation}
As it was mentioned,  the insertion of the shift $\tau$ in the formulas allows to write some relations in the classical  form:
\begin{prop}
Let $\lambda=(\lambda_1,\dots \lambda_m)$ be a composition of number $k$ (the order of parts is important).
Let $a_i=\lambda_1+\dots +\lambda_i$, ($i=1,2,\dots m$). Then
\begin{equation}\label{e_lambda}
e_k(u,\tau)=\sum_{\lambda}\frac{(-1)^{k-m}}{a_1a_2\dots a_m}
p^-_{\lambda_1}(u,\tau)\dots p^-_{\lambda_m}(u,\tau),
\end{equation}
\begin{equation}\label{h_lambda}
h_k(u,\tau)=\sum_{\lambda}\frac{1}{a_1a_2\dots a_m}
p^+_{\lambda_1}(u,\tau) \dots p^+_{\lambda_m}(u,\tau),
\end{equation}
where the sums in both equations are taken over all compositions $\lambda$ of the number $k$.
\end{prop}
\begin{remark}
Compare these formulas with  ($2.14^\prime$) in  Chapter 1.2 of \cite{Mcd}.
\end{remark}
\begin{proof} 
We will prove (\ref{e_lambda}),  the arguments for (\ref{h_lambda})  follow the same lines.
The matrix $B^{-}_k$ can be written as a sum of terms of the form
$$
(P_{m-1, m}\dots P_{a_{m-1}-1, a_{m-1}})\dots(P_{a_{1}-1, a_{1}}\dots P_{1, 2}),
$$
with  permutation matrices $P_{k,l}$, defined by  (\ref{perm}). Each term in this sum corresponds to a decomposition $\lambda$ of number $k$, and the coefficients of these terms in the sum are exactly 
${(-1)^{k-m}}{(a_1a_2\dots a_m)^{-1}}$.
Then from (\ref{e_B}), the elementary  symmetric functions are the sums of 
the products of  terms of the following form:
\begin{equation}\label{PT}
\text{tr}\,
\left(P_{a_i-1,a_i} \dots P_{a_{i-1}-1,a_{i-1}}
T_{a_{i-1}}(u-a_{i-1}+1)  \dots T_{a_{i}}(u-a_{i}+1)\,\right).
\end{equation}
 The following statement can be  checked directly.
 \begin{lemma}\label{RP}
For  any   $k$  matrices $X{(1)},\dots, X{(k)}$ of the size $n\times n$ with the entries in  an associative non-commutative algebra $A$, one has
\begin{equation}\label{P} 
tr\,\left(\,P_{k-1,k}P_{k-2,k-1} \dots P_{1,2}\,(X{(1)})_1\,( X{(2)})_2 \dots(X{(k)})_k\right)=
tr\,(X{(1)} X{(2)}\dots \cdot X{(k)}).
\end{equation}
 \end{lemma}
From Lemma \ref{RP}, the expression in (\ref{PT}) is nothing else but 
$p^{-}_{\lambda_i}(u-a_{i-1}+1)$. Thus, $e_k(u)$ is the sum of terms of the form
$$
{(-1)^{k-m}}{(a_1a_2\dots a_m)^{-1}}p^-_{\lambda_1}(u)
p^-_{\lambda_2}(u-a_1)\dots p^-_{\lambda_m}(u-a_{m-1}),$$
and (\ref{e_lambda}) follows.
\end{proof}

The  following Newton identities   and some of their corollaries are discussed in  \cite {Cherv}, using the technics of so-called  Manin matrices. Here we give  an alternative proof, using the RTT   equation for the Yangian. It is  inspired by the paper  \cite{IOP} on Newton's identities  for RTT algebras with R-matrices  that satisfy Hecke type condition.

\begin{prop}
(Newton's formula)
For any $m=1, 2,\dots$
\begin{equation}\label{rel_ep}
\sum_{k=0}^{m-1}(-1)^{m-k-1}e_{k}(u,\tau)p^{-}_{m-k}(u,\tau)=m\,e_{m}(u,\tau),
\end{equation}
\begin{equation}\label{rel_hh_hp}
\sum_{k=0}^{m-1}h_{k}(u,\tau) p^{+}_{m-k}(u,\tau)=m\, h_{m}(u,\tau).
\end{equation}
\end{prop}
\begin{proof}
By (\ref{proof2-4}),
\begin{equation*}
\begin{aligned}
me_m(u)&
 =\text{tr}\,\left(\,R_{m-1,m}\left(\frac{1}{m-1}\right)\,\,A_{m-1}
 T_1(u)\dots T_m(u-m+1)  \right)\\
 \quad&=\text{tr}\,\left(\,A_{m-1}
 T_1(u)\dots T_m(u-m+1)  \right)\\
 \quad & \quad  \quad  \quad - (m-1)\text{tr}\,\left(\,P_{m-1,m}\,A_{m-1}
 T_1(u)\dots T_m(u-m+1)  \right)\\
\quad & =e_{m-1}(u)p_1(u-m+1)\\
\quad & \quad  \quad  \quad- (m-1)\text{tr}\,\left(\,A_{m-1}
 T_1(u)\dots T_m(u-m+1) P_{m-1,m}\, \right).
 \end{aligned}
\end{equation*}
Applying the  cyclic property of the trace, and the Proposition \ref{AandS}, (c) and (e) to the second term in the last expression, we obtain that
\begin{equation*}
\begin{aligned}
me_m(u)&
 =e_{m-1}(u)\,p_1(u-m+1)\\
\quad & \quad  \quad  \quad -\,\text{tr}\,\left(\,A_{m-2} T_1(u)\dots T_m(u-m+1) P_{m-1,m}\, \right)\\
 \quad & \quad  \quad  \quad+(m-2)\,\text{tr}\,\left(\,A_{m-2} T_1(u)\dots T_m(u-m+1) P_{m-1,m}P_{m-2,m-1}\, \right),
 \end{aligned}
\end{equation*}
and by induction,
\begin{equation}
\begin{aligned}
me_m(u)&
 =e_{m-1}(u)p_1(u-m+1) \\
\quad & - \text{tr}\,\left(\,A_{m-2} T_1(u)\dots T_m(u-m+1) P_{m-1,m}\, \right)\\
 \quad & +\text{tr}\,\left(\,A_{m-3} T_1(u)\dots T_m(u-m+1) P_{m-1,m}P_{m-2,m-1}\, \right) 
 +\dots \\
 \quad &+(-1)^{m-1} \text{tr}\,\left(T_1(u)\dots T_m(u-m+1) P_{m-1,m}\dots P_{1,2}\, \right).
 \end{aligned}
\end{equation}
Applying Lemma \ref{RP} to  the terms of the sum, we  conclude that  each of them has the form  $$(-1)^{m-k-1}e_k(u)p^-_{m-k}(u-k),$$ and the Newton's formula for elementary symmetric functions $e_m(u,\tau)$ follows.

The proof  for homogeneous functions is similar.
\end{proof}

\begin{cor}
 (a)  Coefficients of $\{p^-_k(u)\}$ belong to the Bethe subalgebra $B(n)$. Therefore, they commute. 
 
 (b) 
 \begin{equation*}
m!\,e_m(u)=\det
\begin{pmatrix}
p^-_1(u)\ &  1 & 0 & \dots & 0\\
p^-_2(u) &  p^-_1(u-1)& 2 & \dots & 0\\
\dots &\dots &\dots &\dots &\dots \\
p^-_m(u)&  p^-_{m-1}(u-1)&p^-_{m-2}(u-2)& \dots & p^-_1(u-m+1)\\
\end{pmatrix},
\end{equation*}

\begin{equation*}
m!\,h_m(u)=\det
\begin{pmatrix}
p^+_1(u)&  -1& 0 & \dots & 0\\
p^+_2(u) &  p^+_1(u+1) & -2 & \dots & 0\\
\dots &\dots &\dots &\dots &\dots \\
p^+_m(u)&  p^+_{m-1}(u+1) &p^+_{m-2}(u+2) & \dots & p^+_1(u+m-1)\\
\end{pmatrix},
\end{equation*}

\begin{equation*}
p^-_m(u)=\det
\begin{pmatrix}
e_1(u)&  1& 0 & \dots & 0\\
2\,e_2(u) & e_1(u-1) & 1 & \dots & 0\\
\dots &\dots &\dots &\dots &\dots \\
 m\,e_m(u) & e_{m-1}(u-1)&e_{m-2}(u-2) & \dots & e_1(u-m+1)\\
\end{pmatrix},
\end{equation*}

\begin{equation*}
(-1)^{m-1}\,p^+_m(u)=\det
\begin{pmatrix}
h_1(u)& 1 & 0 & \dots & 0\\
2\,h_2(u) & h_1(u+1)  & 1 & \dots & 0\\
\dots &\dots &\dots &\dots &\dots \\
m\,h_m(u) &  h_{m-1}(u+1)&h_{m-2}(u+2)& \dots & h_1(u+m-1)\\
\end{pmatrix}.
\end{equation*}

\end{cor}

\subsection*{Inverse of the universal differential operator}
 
 Consider the  universal differential operator for XXX model: the formal polynomial in variable $\tau^{-1}$, which is the generating function of the elements $e_k(u)$  (see e.g. \cite{MTV1}, \cite{Tal}):
$$
E(u,\tau)=\sum_{k=0}^{n}(-1)^{k}e_k(u,\tau).
$$
 Using the Newton's identities, it is easy to describe the inverse of this operator.

Namely, define $h^{-}_m(u)$ and $h^{-}_m(u,\tau)$ by the following formulas:
 \begin{equation}\notag
h^-_m(u,\tau):=\tau^{-m}h^-_m(u),
\end{equation}
where 
 \begin{equation}\notag
 \begin{aligned}
m!\, h^{-}_m(u)=
\text{det}\,
\begin{pmatrix}
p^-_1(u) &  -1 &  \dots & 0\\
p^-_2(u+1) &  p^-_1(u+1)  & \dots & 0\\
\dots &\dots &\dots  &\dots \\
p^-_{m-1}(u+m-2)&  p^-_{m-2}(u+m-2)& \dots & -m+1\\
p^-_m(u+m-1)&  p^-_{m-1}(u+m-1)& \dots & p_1^-(u+m-1)\\
\end{pmatrix}.
\end{aligned}
\end{equation}
Let  $$
 H^-(u,\tau)=\sum_{l=0}^{\infty}h^-_l(u,\tau).
 $$
 The following proposition   follows directly from   Newton's identities.
 \begin{prop}
 (a)  The generating functions $H(u,\tau)$, $E(u,\tau)$  satisfy the following identity:
$$
 E(u,\tau)H^-(u+1,\tau)=1
 $$
 (b) The coefficients of the elements $\{ h^{-}_{k}(u)\}$ belong to Bethe subalgebra and commute. The relation to elementary symmetric functions is given by
 $$
 e_k(u)=\text{det}\,(h^{-}_{j-i+1}(u-j+1)).
 $$
  \end{prop}
 One can go further and introduce combinatorial analogues of Schur fucntions:
 \begin{defn} Let $\lambda=(\lambda_1,\dots \lambda_k)$ be a partition of number $m$ of length  at most $k$.
 The {\it Schur funtion }  $s_\lambda(u)$ is the formal series in $u^{-1}$ with coefficients in $Y(n)$, defined by 
\begin{equation}\label{s_def}
 s_\lambda(u):= det\,[h^-_{\lambda_i-i+j}(u-j+1)]_{1\le i,j\le k}.
\end{equation} 
 \begin{prop}
 Let $\lambda^\prime$ be the conjugate partition  to $\lambda$, and assume that it has length at most  $k^{\,\prime}$. Then 
  \begin{equation}\label{s_e}
 s_\lambda(u):= det\,[e_{\lambda^{\prime}_i-i+j}(u)]_{1\le i,j\le k^{\prime}}.
 \end{equation} 

\end{prop}
 \end{defn}
 \begin{proof}
The proof is the same as in  classical case (see \cite{Mcd}, (2.9), ($2.9^\prime$), (3.4), (3.5)).
For any positive number $N$ consider the matrices 
$$
{H^-}=[\,h^-_{i-j}(u-j+1)\,]_{0\le i,j\le N},\quad \quad
{E}=[\,(-1)^{i-j}e_{i-j}(u)\,]_{0\le i,j\le N}.
$$
Here $h^-_{k}(u)=e_k(u)=0$ for any  $k<0$.
The Newton's identities  show that these matrices are inverses of each other. Therefore, each minor of $H^-$ is equal to the complementary cofactor of the transpose of $E$, which implies   the equality of determinants in (\ref{s_def}) and (\ref{s_e}) (c.f.   \cite{Mcd}, formulas
 (2.9), ($2.9^\prime$)).
 \end{proof}

\section* {Connection  to Capelli polynomials and  Shifted Schur functions} \label{Ug}
In this section we show that the proved above identities immediately imply similar relations between Capelli  polynomials and  shifted Schur functions.
The theory  of higher Capelli polynomials is contained in  \cite{Naz}. The detailed account on  shifted symmetric fucntions and their applications is  developed in 
\cite{OO1}. Here we briefly remind  the main definitions, following these two references.

Let
$E=\{e_{ij}\}$ be the  matrix of generators of $\gl_n(\CC)$.
Let $\lambda=(\lambda_1\dots \lambda_k)$ be a partition of number $m$, let $\{c_i\}$ be the set of contents of a column tableau of shape $\lambda$ (see \cite{Naz}  for more details). Consider the  Schur projector  $F_\lambda$   in the tensor power  $(\CC^n)^{\otimes m}$ to  the irreducible 
$\gl_n(\CC)$-component $V_\lambda$.
\begin{defn}The higher Capelly polynomial $c_\lambda(u)$ is a polynomial in  variable $u$ and coefficients in the universal enveloping algebra $U(\gl_n(\CC))$, defined by
\begin{equation}
c_\lambda(u)=\text{tr}(F_\lambda\otimes 1\,(u-c_1+E)_1\dots (u-c_k+E)_k).
\end{equation}
\end{defn}
The coefficients of Capelly polynomials $c_\lambda(u)$  are in  the center of $U(\gl_n(\CC))$. 
The Capelli element  $c_\lambda(u)$  acts in the irreducible representation $V_\mu$ with the highest  weight $\mu$  by multiplication by a scalar, which  is the shifted symmetric function
$s^*_\lambda(\mu+u)$. The  constant coefficients  $\{c_\lambda(0)\}$ form a linear basis of the center of $U(\gl_n(\CC))$. 
In particular, we consider the shifted elementary  functions $e^*_k(u)=s^*_{(1^k)}(\mu+u)$ and  shifted homogeneous symmetric functions $h^*_k(u)=s^*_{(k)}(\mu+u)$, which take the form 
$$
e^*_k(u)=\sum_{1\le i_1<i_2<\dots <i_k<\infty}(\mu_{i_1}+u+k-1)
(\mu_{i_2}+u+k-2)\dots(\mu_{i_k}+u),
$$
$$
h^*_k(u)=\sum_{1\le i_1\le i_2\le \dots \le i_k<\infty}(\mu_{i_1}+u-k+1)
(\mu_{i_2}+u-k+2)\dots(\mu_{i_k}+u).
$$
We  identify the corresponding Capelli elements with their shifted Schur functions, and use the notations $e^*_k(u)$, $h^*_k(u)$
for $c^*_{(1^k)}(u)$ and $c_{(k)}(u)$ respectively.

Let $ev:Y(n)\to U(\gl_n(\CC)) $  be the evaluation homomorphism:
$$
ev: T(u) \mapsto  1+\frac{E}{u}.
$$
Under this map the defined above symmetric functions in $Y(n)$ map to the following Capelli elements:
$$
ev\,(e_k(u))=\frac{e^*_{k}\,{(u-k+1)}}{(u\downarrow k)},\quad 
ev\,(h_k(u))=\frac{h^*_{k}\,{(u+k-1)}}{(u\uparrow k)},
$$
where $$(u\downarrow k)=u(u-1)\dots (u-k+1)\quad \text{ and }
 \quad (u\uparrow k)=u(u+1)\dots (u+k-1).$$

Moreover,
set $$p_m(u)=\text{tr}\,( (E+u)\dots (E+u+m-1)).
$$
Then
$$
ev\,(\,p^-_m(u+m-1)\,)=ev\,(p^+_m(u))=\frac{p_m(u)}{(u\uparrow m)},
$$
and this implies $$
ev\,(h^-_m(u))=ev\,(h_m(u)).$$

The eigenvalue of the central polynomial $p_k(u)\in U(\gl_n(\CC))[u]$ in the irreducible representation $V_\mu$ can be easily found, using the classical formula for the eigenvalues of Casimir operators from \cite{Per-Pop}. The eigenvalue of $\text{tr} \,E^k$ is given by the formula
\begin{equation}
\text{tr} \,E^k(\mu)=\sum_{i=1}^{n}\gamma_i\,m_i^k,
\end{equation}
where 
$$
m_i=\mu_i+n-i,\quad
 \gamma_i=\prod_{j\ne i}\left(
 1-\frac{1}{m_i-m_j}
\right).$$
Accordingly, the shifted symmetric function $p^*_k(u)$ which gives the eigenvalue of $p_k(u)$ is 
$$
p^*_k(u)=\sum_{i=1}^{n}\gamma_i\,(m_i+u)(m_i+u+1)
\dots (m_i+u+k-1).$$

The   combinatorial identities in the Yangian imply immediately  the corresponding relations between Capelli polynomials. Some of them are listed below.

\begin{prop}
\begin{equation}e^*_k(u-k)=\sum_{\lambda} \frac{(-1)^{k-m}}{a_1a_2\dots a_m} p^*_{\lambda_1}(u-a_1)\dots p^*_{\lambda_m}(u-a_m),
\end{equation}

\begin{equation}
h^*_k(u+k-1)=\sum_{\lambda}\frac {1}{a_1a_2\dots a_m} p^*_{\lambda_1}(u)\dots p^*_{\lambda_m}(u+a_{m-1}),
\end{equation}
\begin{equation}\label{eh*}
\sum_{k=0}^{m}(-1)^k e^*_k(u-k+1)\,h^*_{m-k}(u-k)=
\delta_{m,0}.
\end{equation}
\end{prop}
\begin{remark}
 The defined here  functions   $p^*_k(u)$ do not coincide with the shifted power sums under the same notation in \cite {OO1}.

The identity (\ref{eh*}) is similar to the relation (12.18), \cite{OO1} on the generating functions of the elements $e^*_k$, $h^*_k$.
\end{remark}

  \end{document}